\documentclass[11pt]{amsproc}
\usepackage[utf8]{inputenc}
\usepackage{amsthm,amsmath,amssymb}
\usepackage[margin=0.9in]{geometry}
\usepackage{microtype}
\usepackage{tikz}
\usetikzlibrary{arrows.meta,decorations.markings,calc}
\usepackage{float}
\usepackage{xcolor}

\newcommand\dd{\,d}
\newcommand\ii{\mathrm i}
\newtheorem{Lemma}{Lemma}
\newtheorem{MainThm}{Theorem}


\title[Fourier harmonics and pentagonal rigidity]{Fourier Harmonics and Pentagonal Rigidity in the Planar Ulam Floating Body Problem}

\author{Oleg Asipchuk}
\address{Oleg Asipchuk, University of Cincinnati, Department of Mathematical Sciences, Cincinnati, OH 45221, USA}
\email{asipchah@ucmail.uc.edu}

\author{Maksim Kosmakov}
\address{Maksim Kosmakov, University of Toledo, Department of Mathematics and Statistics, Toledo, OH 43606, USA}
\email{Maksim.Kosmakov@UToledo.Edu}

\begin{document}

\begin{abstract}
We give a new proof that the disk is the only strictly convex planar body floating
in equilibrium in every orientation at perimetral densities
$\sigma=1/5$ and $\sigma=2/5$.
The proof uses the Fourier harmonics of the chord-angle function
and the closure relations of the associated pentagonal carousel.
Their combination gives inequalities that force the body to be a disk.
\end{abstract}
\maketitle

\section{Introduction}

Ulam's Problem~19 asks whether a homogeneous body that floats
in equilibrium in every orientation must be a Euclidean sphere
\cite{ScottishBook}. Auerbach studied the planar analogue
and showed that noncircular floating bodies exist at physical
density $1/2$, with boundaries given by the classical Zindler
curves~\cite{Auerbach,Zindler}. More recently, Ryabogin gave a
negative answer to the original problem in three dimensions
\cite{Ryabogin2022}; see also the survey \cite{ARSY}. Furthermore, the result of Auerbach was further extended in \cite{MR2026} to include certain fractal curves.

In the planar setting, Auerbach's approach relates the
flotation condition to a family of chords of constant length.
These chords divide the boundary into arcs of lengths
$\sigma P$ and $(1-\sigma)P$, where $P$ is the perimeter
and $\sigma\in(0,1/2]$ is called the \emph{perimetral density}.
Thus, the planar problem can be studied through the geometry
of equal chords on convex curves.

When $\sigma=k/n$ is rational, with $k$ and $n$ coprime,
successive advances by arclength $kP/n$ trace an equilateral
$n$-gon. As the starting point moves along the boundary,
this polygon moves with it, forming a \emph{Zindler carousel}.
Bracho, Montejano, and Oliveros developed this approach
and established rigidity for $\sigma=1/5$ and $\sigma=2/5$
without assuming convexity \cite{BMO1,BMO2}.
Their proof studies the dynamics of equilateral pentagons and uses estimates for the associated period function, partially supported by numerical evidence.

Earlier characterizations of circles by their chords
include the work of Salkowski \cite{Salkowski1934}.
Bracho, Montejano, and Oliveros also established rigidity
for $\sigma=1/3$ and $\sigma=1/4$ \cite{BMO2}.
More recently, rigidity was proved for $\sigma=1/6$
and, under central symmetry, for $\sigma=1/8$ and
$\sigma=3/8$ \cite{AKZ16,AKZ18}.
 On the other hand, Wegner constructed families of noncircular floating bodies
for certain physical densities different from $1/2$
\cite{Wegner2003,Wegner2007}. Here the physical density
is distinct from the perimetral density $\sigma$ considered
in the present paper.

In this paper, we give a short Fourier-analytic proof for
the pentagonal densities $\sigma=1/5$ and $\sigma=2/5$.
Instead of studying the period of the carousel flow, we use
two Fourier decompositions. The first, along the boundary,
relates the Fourier coefficients of the chord angle to those
of its sine. The second, over the five edges of the carousel,
uses the fact that those edges form a closed polygon.
Together, these relations force the chord angle to be
constant, and thus the boundary to be a circle.

Our main result is the following.

\begin{MainThm}\label{thm:main}
Let $K\subset\mathbb R^2$ be a strictly convex planar body
with $C^1$ boundary. If $K$ floats in equilibrium in every
orientation with perimetral density $\sigma=1/5$ or
$\sigma=2/5$, then $K$ is a disk.
\end{MainThm}

The paper is organized as follows. Section~2 introduces the carousel notation and Fourier
identities. Section~\ref{sec:angle-bounds} establishes the
geometric angle bounds, and
Section~\ref{sec:fourier-proof} completes the proof.

\medskip
\noindent
{\it Acknowledgments}.
The authors would like to thank Dmitry Ryabogin for bringing this question to their attention and Pavel Zatitskii for many insightful mathematical discussions. 

\medskip
\noindent
{\it AI Declaration}. During the preparation of this manuscript, the authors used ChatGPT, developed by OpenAI,  as an assistive tool for editing and typesetting the manuscript and for checking the correctness of mathematical arguments and proofs. All mathematical results, however, were obtained independently by the authors without the assistance of artificial intelligence.  

\section{The carousel and its Fourier identities}

Let $K$ satisfy the hypotheses of Theorem~\ref{thm:main}, let
$P=|\partial K|$, and let
$\gamma:\mathbb R/P\mathbb Z\to\partial K$ be a positively oriented arc-length
parametrization.  Write $\sigma=k/5$ with $k\in\{1,2\}$ and set $h=kP/5$.
Auerbach's criterion gives a common length for the flotation chords.  After a
similarity, normalize this length to $2$.  The corresponding unit chord vector
is therefore
\begin{equation}\label{eq:edge-function}
 e(t):=\frac{\gamma(t+h)-\gamma(t)}2,
 \qquad |e(t)|=1.
\end{equation}
Strict convexity places the chord direction between the tangent directions at
its endpoints.  Since the chord length is constant, it makes the same angle
$\alpha(t)\in(0,\pi)$ with the two endpoint tangents.  Thus, under identification of $\mathbb R^2$ with $\mathbb C$, we have 
\begin{equation}\label{eq:endpoint-vectors}
 \gamma'(t)=e(t)e^{-\ii\alpha(t)},
 \qquad
 \gamma'(t+h)=e(t)e^{\ii\alpha(t)}.
\end{equation}
Differentiating \eqref{eq:edge-function} and using
\eqref{eq:endpoint-vectors} gives
\begin{equation}\label{eq:edge-derivative}
 e'(t)=\ii\sin\alpha(t)\,e(t).
\end{equation}
Choose a continuous lift $\theta:\mathbb R\to\mathbb R$ of the argument of $e$, so that
$e(t)=e^{\ii\theta(t)}$.
Then
\begin{equation}\label{eq:theta-derivative}
 \theta'(t)=\sin\alpha(t)>0.
\end{equation}
Since $\gamma'$ and $e$ are $P$-periodic, \eqref{eq:endpoint-vectors} and
$0<\alpha<\pi$ imply that $\alpha$ is $P$-periodic.  The tangent direction
$\theta-\alpha$ winds once around the positively oriented strictly convex
curve, and thus
\begin{equation}\label{eq:theta-winding}
 \theta(t+P)=\theta(t)+2\pi.
\end{equation}

Evaluating the common tangent at $\gamma(t+h)$ in the two endpoint formulas gives
\[e(t)e^{\ii\alpha(t)}=e(t+h)e^{-\ii\alpha(t+h)}.\]
Thus, the turn from $e(t)$ to $e(t+h)$ is $\alpha(t)+\alpha(t+h)$.  This is the
actual positive turn of two successive carousel edges, so strict convexity
gives
\begin{equation}\label{eq:theta-shift}
 \theta(t+h)-\theta(t)=\alpha(t)+\alpha(t+h),
 \qquad 0<\alpha(t)+\alpha(t+h)<\pi.
\end{equation}
The five shifted copies of \eqref{eq:theta-shift}
form a cyclic linear system for the five values
$\alpha(t+jh)$. Since the number of equations is odd,
the system is invertible. As $\theta\in C^1$,
it follows that $\alpha\in C^1$.  Differentiating and using
\eqref{eq:theta-derivative} gives
\begin{equation}\label{eq:alpha-delay}
 \alpha'(t)+\alpha'(t+h)
 =\sin\alpha(t+h)-\sin\alpha(t).
\end{equation}

The five successive chords form the Zindler carousel.  Set
\begin{equation}\label{eq:carousel-variables}
 \begin{aligned}
 v_j(t)&:=\gamma(t+(j-1)h),&
 e_j(t)&:=e(t+(j-1)h)=\frac{v_{j+1}-v_j}{2},\\
 \theta_j(t)&:=\theta(t+(j-1)h),&\alpha_j(t)&:=\alpha(t+(j-1)h).
 \end{aligned}
\end{equation}
Here $v_j,e_j,\alpha_j$ are read modulo $5$, while $\theta_j$ is the continuous
lift, so $\theta_{j+5}=\theta_j+2\pi k$.
Since $v_6=v_1$, we have
\begin{equation}\label{eq:edge-closure}
 e_1+e_2+e_3+e_4+e_5=0.
\end{equation}

By \eqref{eq:theta-shift}, the positive turn from
$e_{j-1}$ to $e_j$ is
$\alpha_{j-1}+\alpha_j\in(0,\pi)$.
Thus, the angle between the two incident edges
at $v_j$ is
\begin{equation}\label{eq:xalpha}
 x_j:=\pi-\alpha_{j-1}-\alpha_j\in(0,\pi).
\end{equation}
In particular,
\begin{equation}\label{eq:triple-turn}
 \theta_j-\theta_{j-1}=\alpha_{j-1}+\alpha_j=\pi-x_j.
\end{equation}

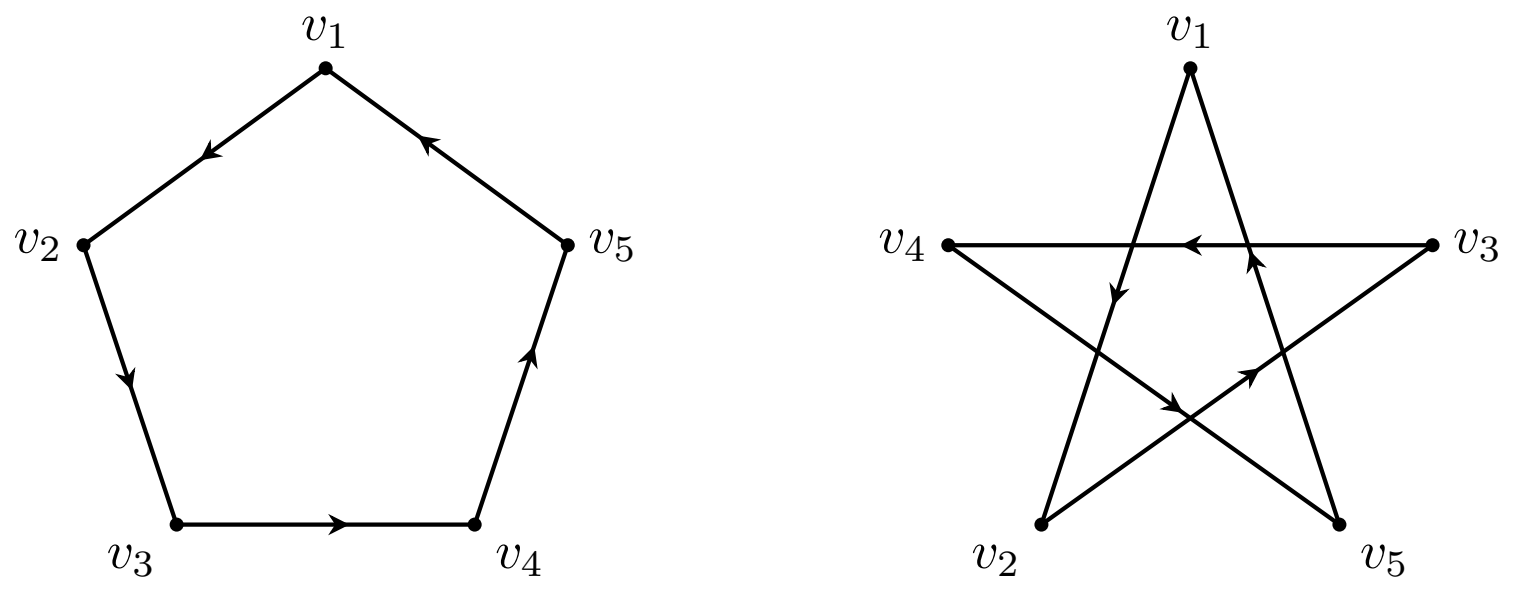
\begin{figure}[t]
\centering
\begin{tikzpicture}[scale=0.88,>=Stealth]
  \begin{scope}[xshift=0cm]
    \coordinate (P1) at (0,1.95);
    \coordinate (P2) at (-1.82,0.62);
    \coordinate (P3) at (-1.12,-1.48);
    \coordinate (P4) at (1.12,-1.48);
    \coordinate (P5) at (1.82,0.62);
    \draw[thick] (P1)--(P2);
    \draw[thick] (P2)--(P3);
    \draw[thick] (P3)--(P4);
    \draw[thick] (P4)--(P5);
    \draw[thick] (P5)--(P1);
    \draw[-{Stealth[length=4pt,width=4pt]}] ($(P1)!0.35!(P2)$)--($(P1)!0.52!(P2)$);
    \draw[-{Stealth[length=4pt,width=4pt]}] ($(P2)!0.35!(P3)$)--($(P2)!0.52!(P3)$);
    \draw[-{Stealth[length=4pt,width=4pt]}] ($(P3)!0.38!(P4)$)--($(P3)!0.58!(P4)$);
    \draw[-{Stealth[length=4pt,width=4pt]}] ($(P4)!0.48!(P5)$)--($(P4)!0.64!(P5)$);
    \draw[-{Stealth[length=4pt,width=4pt]}] ($(P5)!0.45!(P1)$)--($(P5)!0.62!(P1)$);
    \fill (P1) circle (1.5pt) node[font=\small,above] {$v_1$};
    \fill (P2) circle (1.5pt) node[font=\small,left] {$v_2$};
    \fill (P3) circle (1.5pt) node[font=\small,below left] {$v_3$};
    \fill (P4) circle (1.5pt) node[font=\small,below right] {$v_4$};
    \fill (P5) circle (1.5pt) node[font=\small,right] {$v_5$};
    \node at (0,-2.55) {$k=1$: pentagon};
  \end{scope}
  \begin{scope}[xshift=6.5cm]
    \coordinate (Q1) at (0,1.95);
    \coordinate (Q2) at (-1.12,-1.48);
    \coordinate (Q3) at (1.82,0.62);
    \coordinate (Q4) at (-1.82,0.62);
    \coordinate (Q5) at (1.12,-1.48);
    \draw[thick] (Q1)--(Q2);
    \draw[thick] (Q2)--(Q3);
    \draw[thick] (Q3)--(Q4);
    \draw[thick] (Q4)--(Q5);
    \draw[thick] (Q5)--(Q1);
    \draw[-{Stealth[length=4pt,width=4pt]}] ($(Q1)!0.35!(Q2)$)--($(Q1)!0.52!(Q2)$);
    \draw[-{Stealth[length=4pt,width=4pt]}] ($(Q2)!0.38!(Q3)$)--($(Q2)!0.56!(Q3)$);
    \draw[-{Stealth[length=4pt,width=4pt]}] ($(Q3)!0.35!(Q4)$)--($(Q3)!0.52!(Q4)$);
    \draw[-{Stealth[length=4pt,width=4pt]}] ($(Q4)!0.42!(Q5)$)--($(Q4)!0.60!(Q5)$);
    \draw[-{Stealth[length=4pt,width=4pt]}] ($(Q5)!0.42!(Q1)$)--($(Q5)!0.60!(Q1)$);
    \fill (Q1) circle (1.5pt) node[font=\small,above] {$v_1$};
    \fill (Q2) circle (1.5pt) node[font=\small,below left] {$v_2$};
    \fill (Q3) circle (1.5pt) node[font=\small,right] {$v_3$};
    \fill (Q4) circle (1.5pt) node[font=\small,left] {$v_4$};
    \fill (Q5) circle (1.5pt) node[font=\small,below right] {$v_5$};
    \node at (0,-2.55) {$k=2$: index-two pentagon};
  \end{scope}
\end{tikzpicture}
\caption{Pentagonal carousels. Arrows indicate the order
of the equal chords.}
\label{fig:carousels}
\end{figure}

The sum of the five shifted copies of \eqref{eq:theta-shift}, together with
$5h=kP$, is
\[
 2\sum_{j=1}^5\alpha_j
 =\theta(t+kP)-\theta(t)=2\pi k.
\]
Thus,
\begin{equation}\label{eq:sum}
 \sum_{j=1}^5\alpha_j=k\pi,
 \qquad
 \sum_{j=1}^5x_j=(5-2k)\pi.
\end{equation}
Averaging the first identity in \eqref{eq:sum} over one perimeter gives
$\mu:=P^{-1}\int_0^P\alpha(t)\,\dd t=k\pi/5$.

Equation~\eqref{eq:alpha-delay} gives
\begin{equation}\label{eq:angle-ode}
 \alpha_{j-1}'+\alpha_j'
 =\sin\alpha_j-\sin\alpha_{j-1},
 \qquad j\pmod5.
\end{equation}
For a $P$-periodic function set
\[
 \widehat f(q):=\frac1P\int_0^P f(t)e^{-2\pi\ii qt/P}\,\dd t,
 \qquad \langle f, g\rangle:=\frac1P\int_0^P f(t)\overline{g}(t)\,\dd t,\qquad 
 \|f\|_2^2:=\frac1P\int_0^P|f(t)|^2\,\dd t.
\]
With this normalization, Parseval's identity is
\begin{equation}\label{eq:parseval-time}
 \|f\|_2^2=\sum_{q\in\mathbb Z}|\widehat f(q)|^2.
\end{equation}
From \eqref{eq:theta-derivative} and \eqref{eq:theta-winding}, we have
\[
 s:=\widehat{\sin\alpha}(0)
 =\frac1P\int_0^P\sin\alpha(t)\,\dd t
 =\frac{2\pi}{P}.
\]
For $q\ne0$, the comparison of the $q$-th Fourier coefficients of \eqref{eq:alpha-delay} gives
\begin{equation}\label{eq:fourier-raw}
 \frac{2\pi\ii q}{P}\Bigl(1+e^{2\pi\ii kq/5}\Bigr)\widehat\alpha(q)
 =\Bigl(e^{2\pi\ii kq/5}-1\Bigr)\widehat{\sin\alpha}(q).
\end{equation}
Since $2\pi/P=s$, for $q\ne0$ we obtain
\begin{equation}\label{eq:fourier-multiplier}
\begin{cases} \widehat{\sin\alpha}(q)=qs\cot\frac{\pi kq}{5}\,\widehat\alpha(q)\quad&5\nmid q,\\[1mm]
\widehat\alpha(q)=0\quad&5\mid q.
\end{cases}
\end{equation}

We shall also use the five-point Fourier transform.  Put
$\rho=e^{2\pi\ii/5}$.  For $c=(c_1,\ldots,c_5)$ define
\[
 \widetilde c_m:=\frac15\sum_{j=1}^5c_j\rho^{-m(j-1)},
 \qquad m=0,\ldots,4.
\]
Then
\begin{equation}\label{eq:cyclic-parseval}
 c_j=\sum_{m=0}^4\widetilde c_m\rho^{m(j-1)},
 \qquad
 \frac15\sum_{j=1}^5|c_j|^2=\sum_{m=0}^4|\widetilde c_m|^2.
\end{equation}
Apply the five-point transform to the angle vector and write
\[
 z_m(t):=\frac15\sum_{j=1}^5\alpha_j(t)\rho^{-m(j-1)},
 \qquad m=0,\ldots,4.
\]
Then $z_0=\mu$, $z_{5-m}=\overline{z_m}$, and
\begin{equation}\label{eq:angle-cyclic-mass}
 \|\alpha-\mu\|_2^2=2(\|z_1\|_2^2+\|z_2\|_2^2).
\end{equation}
Since $z_m(t+h)=\rho^m z_m(t)$, its temporal Fourier frequencies satisfy
$kq\equiv m\pmod5$.

\begin{Lemma}\label{lem:constant-angle}
If $z_1=z_2=0$ identically, then $\gamma$ is a circle.
\end{Lemma}

\begin{proof}
By \eqref{eq:angle-cyclic-mass}, $\alpha=\mu=k\pi/5$.
With $\kappa:=\sin(k\pi/5)>0$,
\eqref{eq:edge-derivative} gives
$e(t)=e(0)e^{\ii\kappa t}$, while
\eqref{eq:endpoint-vectors} gives
$\gamma'(t)=e^{-\ii k\pi/5}e(t)$.
Thus, $\gamma$ parametrizes a circle.
\end{proof}

\section{Pentagonal angle bounds}\label{sec:angle-bounds}

\subsection{Index-two case}\label{sec:index-two}
\begin{Lemma}\label{lem:k2-angles}
If $k=2$, then
\[
 \frac\pi3<\alpha_j<\frac{2\pi}{3},
 \qquad s>\frac{\sqrt3}{2}.
\]
\end{Lemma}

\begin{proof}
By cyclic relabeling, it is sufficient to consider $\alpha_1$. Since $\sum_j\alpha_j=2\pi$, the relations \eqref{eq:xalpha} give $x_3+x_5=\alpha_1$. Set $A:=v_4-v_2$ and $B:=v_1-v_4$, so $|A+B|=2$.

\begin{figure}[ht]
\centering
\includegraphics[width=.4\textwidth]{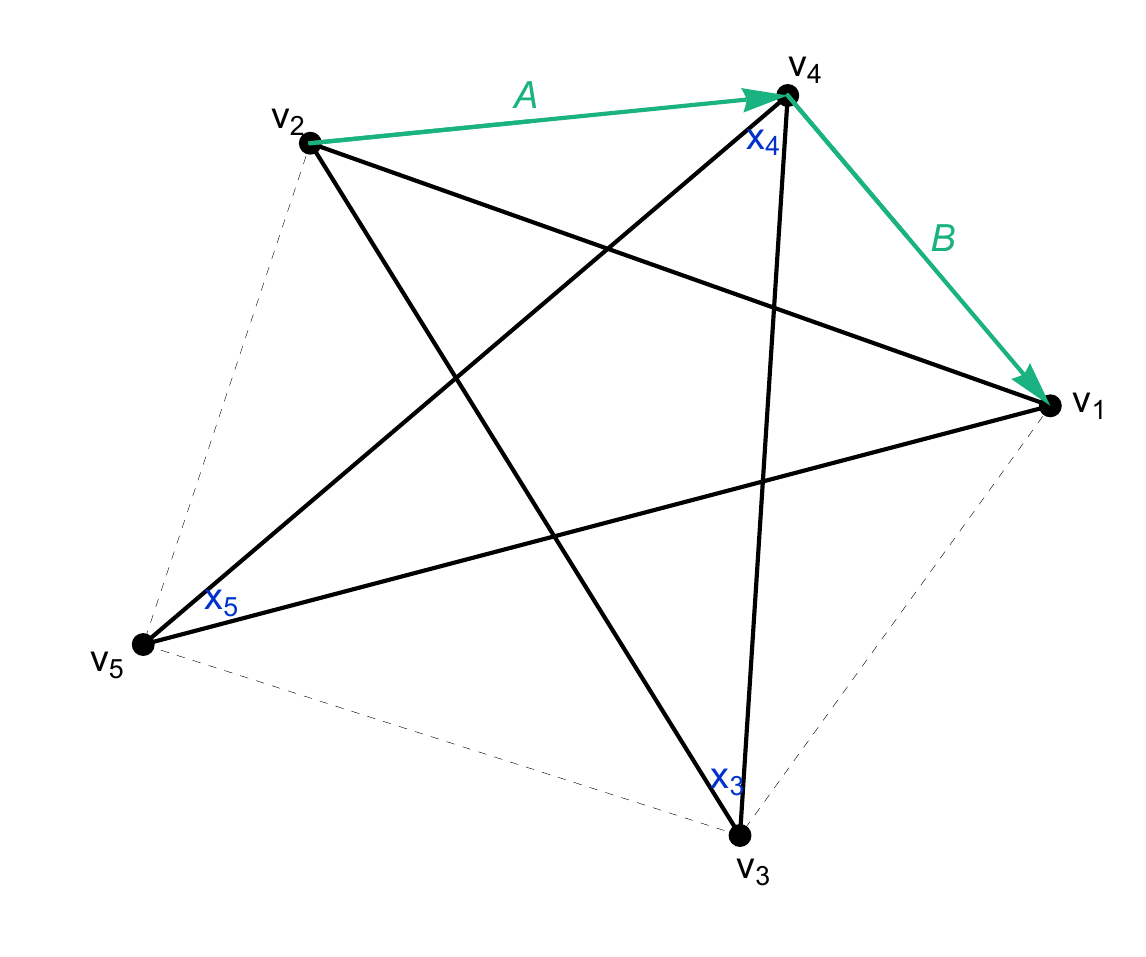}
\caption{The geometry used in Lemma~\ref{lem:k2-angles}.}
\label{fig:lemma2}
\end{figure}

The triangles $v_2v_3v_4$ and $v_4v_5v_1$ are isosceles (see Figure~\ref{fig:lemma2}), with apex angles $x_3$ and $x_5$. Thus,
\begin{equation}\label{eq:sides-A-B}
|A|=4\sin\frac{x_3}{2},\qquad |B|=4\sin\frac{x_5}{2},
\end{equation}
and their base angles are $(\pi-x_3)/2$ and $(\pi-x_5)/2$. Therefore,
\[
\angle v_2v_4v_1=\pi-\left(x_4+\frac{x_3+x_5}{2}\right).
\]
Since $\sum_j x_j=\pi$ and $x_j>0$, we have
\[
 0<\frac{\alpha_1}{2}=\frac{x_3+x_5}{2}<x_4+\frac{x_3+x_5}{2}<\pi.
\]
The law of cosines gives
\[
4=|A|^2+|B|^2+2|A||B|\cos\left(x_4+\frac{x_3+x_5}{2}\right)
 <|A|^2+|B|^2+2|A||B|\cos\frac{x_3+x_5}{2}.
\]
Using  $\sin^2u+\sin^2v+2\sin u\sin v\cos(u+v)=\sin^2(u+v)$ and \eqref{eq:sides-A-B}, we have
\[
4<16\sin^2\frac{x_3+x_5}{2}=16\sin^2\frac{\alpha_1}{2},
\]
so $\alpha_1>\pi/3$. The cyclic relabeling gives the same lower bound for all $j$. Since $\alpha_j+\alpha_{j+1}=\pi-x_{j+1}<\pi$ and $\alpha_{j+1}>\pi/3$, we also have $\alpha_j<2\pi/3$. Thus, $\sin\alpha_j>\sqrt3/2$ for every $j$, which gives $s>\sqrt3/2$.

\end{proof}

\subsection{Index-one case}\label{sec:pentagon}
\begin{Lemma}\label{lem:k1-angles}
If $k=1$, then
\[
 0<\alpha_j<\frac\pi2,
 \qquad
 \frac9{20}<s\le\sin\frac\pi5.
\]
\end{Lemma}

\begin{proof}
By cyclic relabeling, it is sufficient to consider $\alpha_1$.
Rotate the edge vectors so that $e_1=1$.  The successive edge turns give
\[
 e_5=e^{-\ii(\alpha_5+\alpha_1)},\qquad
 e_2=e^{\ii(\alpha_1+\alpha_2)}.
\]
Since $\sum_j\alpha_j=\pi$, the angle between $e_3$ and $e_4$ is
$\alpha_3+\alpha_4=\pi-\alpha_5-\alpha_1-\alpha_2$.
The closure relation therefore gives
\begin{equation}\label{eq:k1-three-edge}
2\sin\frac{\alpha_5+\alpha_1+\alpha_2}{2} =|e_3+e_4|=|e_1+e_2+e_5|=\left|1+e^{-\ii(\alpha_5+\alpha_1)}+e^{\ii(\alpha_1+\alpha_2)}\right|.
\end{equation}

We first show that $\alpha_5+2\alpha_1+\alpha_2<\pi$. Indeed,
\[
\left|1+e^{-\ii(\alpha_5+\alpha_1)}+e^{\ii(\alpha_1+\alpha_2)}\right|^2
=1+8\cos\frac{\alpha_5+\alpha_1}{2}\cos\frac{\alpha_1+\alpha_2}{2}\cos\frac{\alpha_5+2\alpha_1+\alpha_2}{2}.
\]
Since $\alpha_5+\alpha_1<\pi$ and $\alpha_1+\alpha_2<\pi$, the first two cosine factors are positive. If $\alpha_5+2\alpha_1+\alpha_2\ge\pi$, the modulus on the left
is at most $1$.  But $\alpha_5+\alpha_1+\alpha_2>\frac{\alpha_5+2\alpha_1+\alpha_2}{2}\ge\frac{\pi}{2}$, so
\eqref{eq:k1-three-edge} makes the same modulus larger than $\sqrt2$,
a  contradiction. Thus, $\alpha_5+2\alpha_1+\alpha_2<\pi$, and in particular $\alpha_1<\pi/2$.

Since the unit vectors $e_2$ and $e_5$ make the angle
$\alpha_5+2\alpha_1+\alpha_2$, we have
\[
 |e_2+e_5|=2\cos\frac{\alpha_5+2\alpha_1+\alpha_2}{2}
 <2\cos\alpha_1.
\]
Thus, by the triangle inequality and \eqref{eq:k1-three-edge},
\[
2\sin\frac{\alpha_1}{2}<2\sin\frac{\alpha_5+\alpha_1+\alpha_2}{2}=|e_1+e_2+e_5|\le1+2\cos\frac{\alpha_5+2\alpha_1+\alpha_2}{2}<1+2\cos\alpha_1.
\]
The last inequality implies
$\sqrt{2-2\cos\alpha_1}<1+2\cos\alpha_1$, hence
$4\cos^2\alpha_1+6\cos\alpha_1-1>0$.  Consequently,
\[
 \cos\alpha_j>\frac{\sqrt{13}-3}{4}, \qquad j=1,\ldots,5.
\]
Put $A:=\arccos((\sqrt{13}-3)/4)$.  We have $\pi/3<A<\pi/2$.
By concavity of sine, its sum subject to $0\le a_i\le A$ and
$\sum_i a_i=\pi$ is minimized at an extreme point, namely a
permutation of $(A,A,\pi-2A,0,0)$.  Since each $\alpha_i\in(0,A)$,
\[
 \sum_{i=1}^5\sin\alpha_i
 >2\sin A+\sin(\pi-2A)
 =2\sin A(1+\cos A)>\frac94.
\]
Here, the last inequality is a direct evaluation at
$\cos A=(\sqrt{13}-3)/4$.  Averaging over the five shifts gives
$s>9/20$.

Finally, $P^{-1}\int_0^P\alpha(t)\,\dd t=\pi/5$.  The concavity of
sine on $[0,\pi]$ and Jensen's inequality give $s\le\sin(\pi/5)$.
\end{proof}

\section{Fourier comparison and rigidity}\label{sec:fourier-proof}

\begin{Lemma}\label{lem:common-closure}
 Let $r_m:=\cot(m\pi/5)|z_m|$, $m=1,2$. Then for $k=1,2$,
\begin{equation}\label{eq:common-closure}
 r_k\le r_1^2+r_2^2.
\end{equation}
\end{Lemma}

\begin{proof}
We now use closure equation~\eqref{eq:edge-closure}.  Introduce the angular deviation $\eta_j$ of $e_j$ from the
corresponding edge of the   regular pentagon of index $k$ by
\[
\eta_j:=\theta_j-\frac{2\pi k}{5}(j-1).
\]
Let $\overline\eta:=\frac15\sum_{j=1}^5\eta_j$ and $\delta_j:=\eta_j-\overline\eta$, then
\[
\theta_j=\frac{2\pi k}{5}(j-1)+\overline\eta+\delta_j, \qquad \sum_{j=1}^5\delta_j=0,
\]
and therefore
\[
e_j=e^{\ii\overline\eta}\rho^{k(j-1)}e^{\ii\delta_j}.
\]
 By
\eqref{eq:triple-turn},

\begin{equation}\label{eq:common-turn-difference}
 \delta_j-\delta_{j-1}=\alpha_{j-1}+\alpha_j-\frac{2\pi k}{5}.
\end{equation}
 For the  five-point Fourier expansions
\[
\delta_j=\sum_{m=0}^4\widetilde\delta_m\rho^{m(j-1)},
\qquad \alpha_j=\sum_{m=0}^4 z_m\rho^{m(j-1)},
\]
\[
\delta_{j-1} =\sum_{m=0}^4\rho^{-m}\widetilde\delta_m\rho^{m(j-1)}, \qquad
\alpha_{j-1} =\sum_{m=0}^4\rho^{-m}z_m\rho^{m(j-1)},
\]
  substitution  into \eqref{eq:common-turn-difference}
gives
\[
 (1-\rho^{-m})\widetilde\delta_m=(1+\rho^{-m})z_m,
 \qquad m=1,\ldots,4.
\]
Using $|(1+\rho^{-m})/(1-\rho^{-m})|=\cot(m\pi/5)$ for $m=1,2$,
we obtain $|\widetilde\delta_m|=r_m$. 
Since the $\delta_j$ are real, we have
$\widetilde\delta_1=\overline{\widetilde\delta_4}$ and 
 $\widetilde\delta_2=\overline{\widetilde\delta_3}$,
while $\widetilde\delta_0=0$, because $\sum_j\delta_j=0$. Thus, the finite Parseval's identity gives
\begin{equation}\label{eq:common-edge-parseval}
 \frac12\sum_{j=1}^5\delta_j^2=5(r_1^2+r_2^2).
\end{equation}

By closure, $\sum_j\rho^{k(j-1)}e^{\ii\delta_j}=0$, while
$\sum_j\rho^{k(j-1)}=0$.  Subtracting these identities and isolating the linear term gives
\[
 5\ii\widetilde\delta_{5-k}=\ii\sum_j\rho^{k(j-1)}\delta_j
 =-\sum_j\rho^{k(j-1)}(e^{\ii\delta_j}-1-\ii\delta_j).
\]
Since the sum on the left-hand side has modulus $5r_k$
and $|e^{\ii x}-1-\ii x|\le x^2/2$, \eqref{eq:common-edge-parseval} gives
\[
 5r_k\le\frac12\sum_{j=1}^5\delta_j^2=5(r_1^2+r_2^2),
\]
which proves \eqref{eq:common-closure}.
\end{proof}

\begin{proof}[Proof of Theorem~\ref{thm:main}]
\medskip
\textbf{Case $k=2$.} We have $\mu=2\pi/5$ and  Lemma~\ref{lem:k2-angles} gives $\alpha\in(\pi/3,2\pi/3)$.  Since $|\cos x|\le1/2$ on this interval,
\begin{equation}\label{eq:k2-sine-contraction}
 \|\sin\alpha-s\|_2\le\|\sin\alpha-\sin\mu\|_2
 \le\frac12\|\alpha-\mu\|_2.
\end{equation}

Since $h=2P/5$, the identity $z_m(t+h)=\rho^mz_m(t)$ implies that
$\widehat{z_m}(q)\ne0$ only if $2q\equiv m\pmod5$.  Thus, the temporal
frequencies of $z_1$ satisfy $q\equiv3\pmod5$, while those of
$z_4=\overline{z_1}$ satisfy $q\equiv2\pmod5$.  On these residue
classes, the multiplier in \eqref{eq:fourier-multiplier} satisfies
\[
 \left|qs\cot\frac{2\pi q}{5}\right|\ge2s\cot\frac\pi5>1,
\]
where the last inequality follows from $s>\sqrt3/2$.
Therefore, \eqref{eq:fourier-multiplier} and Parseval's identity give
\begin{equation}\label{eq:k2-multiplier-bound}
 \|\sin\alpha-s\|_2^2
 =\sum_{q\ne0}|\widehat{\sin\alpha}(q)|^2
 \ge\sum_{q\equiv2,3\;(\mathrm{mod}\,5)}|\widehat{\sin\alpha}(q)|^2 \ge\sum_{q\equiv2,3\;(\mathrm{mod}\,5)}|\widehat\alpha(q)|^2
 =\|z_1\|_2^2+\|z_4\|_2^2=2\|z_1\|_2^2.
\end{equation}
On the other hand, \eqref{eq:k2-sine-contraction} and
\eqref{eq:angle-cyclic-mass} give
\[
 \|\sin\alpha-s\|_2^2
 \le\frac14\|\alpha-\mu\|_2^2
 =\frac12(\|z_1\|_2^2+\|z_2\|_2^2).
\]
Comparing this with \eqref{eq:k2-multiplier-bound} gives
\begin{equation}\label{eq:k2-mass-separation}
 \|z_2\|_2^2\ge3\|z_1\|_2^2.
\end{equation}
Since
$\alpha_j>\pi/3$ and $\sum_j\alpha_j=2\pi$, we may write
\[
 z_2=\frac15\sum_{j=1}^5\left(\alpha_j-\frac\pi3\right)
 \rho^{-2(j-1)},\quad\text{so}\quad |z_2|\le\frac\pi{15}.
\]
It follows that $r_2=\cot(2\pi/5)|z_2|<1/10$.
For $k=2$, Lemma~\ref{lem:common-closure} therefore gives
\[
 r_2(1-r_2)\le r_1^2, \quad\text{so}\quad 9r_2^2\le r_1^2.
\]
Since $r_m=\cot(m\pi/5)|z_m|$ and
$\cot^2(\pi/5)/\cot^2(2\pi/5)=9+4\sqrt5<18$, we obtain
\begin{equation*} 
 |z_2|^2\le2|z_1|^2,
 \quad\text{so}\quad \|z_2\|_2^2\le2\|z_1\|_2^2.
\end{equation*}
Thus,  both cyclic modes vanish, and Lemma~\ref{lem:constant-angle} completes the proof.

\medskip
\textbf{Case $k=1$.}
We have $\mu=\pi/5$ and  by
Lemma~\ref{lem:k1-angles}, $0<\alpha_j<\pi/2$. Using
$\cot^2(\pi/5)+\cot^2(2\pi/5)=2$, Cauchy--Schwarz, and
\eqref{eq:cyclic-parseval}, we obtain
\[
(r_1+r_2)^2 =\left(\cot\frac{\pi}{5}|z_1| +\cot\frac{2\pi}{5}|z_2|\right)^2\le 2(|z_1|^2+|z_2|^2)
 =\frac15\sum_{j=1}^5(\alpha_j-\mu)^2
 <\frac{3\pi^2}{50}<1,
\]
where we used $\alpha_j^2<(\pi/2)\alpha_j$ and
$\sum_j\alpha_j=\pi$. If $r_1>r_2$, the closure estimate \eqref{eq:common-closure} would give
$ r_1\le r_1^2+r_2^2<r_1(r_1+r_2)<r_1$, thus necessarily $r_1\le r_2$. Since
$\tan(\pi/5)/\tan(2\pi/5)=\sqrt5-2$, this gives
\begin{equation}\label{eq:k1-geometric-bound}
 |z_1|\le(\sqrt5-2)|z_2|,
 \quad \text{so}\quad\|z_1\|_2^2\le(\sqrt5-2)^2\|z_2\|_2^2.
\end{equation}
For $0<\alpha<\pi/2$, concavity of sine gives
\[
 \frac25<\frac{\sin\alpha-\sin\mu}{\alpha-\mu}<1,
\]
where the quotient is interpreted as $\cos\mu$ at $\alpha=\mu$. Therefore,
\[
 \bigl(\sin\alpha-\sin\mu-\tfrac25(\alpha-\mu)\bigr)
 \bigl(\sin\alpha-\sin\mu-(\alpha-\mu)\bigr)\le0.
\]
Integration over one period gives
\[
0\ge\|\sin\alpha-\sin\mu\|_2^2-\frac75\langle\sin\alpha-\sin\mu,\alpha-\mu\rangle+\frac25\|\alpha-\mu\|_2^2,
\]
and therefore Parseval's identity gives
\[
0\ge(s-\sin\mu)^2+\sum_{q\ne0}\left(|\widehat{\sin\alpha}(q)|^2-\frac75\operatorname{Re}\bigl( \widehat{\sin\alpha}(q)\overline{\widehat\alpha(q)}\bigr)+\frac25|\widehat\alpha(q)|^2 \right).
\]
By \eqref{eq:fourier-multiplier}, it can be rewritten as 
\[
0\ge{}(s-\sin\mu)^2
 +\sum_{\substack{q\ne0\\5\mid q}}  |\widehat{\sin\alpha}(q)|^2+
 \sum_{5\nmid q}\left(qs\cot\dfrac{\pi q}{5} -\frac25\right)\left(qs\cot\dfrac{\pi q}{5}-1\right)|\widehat\alpha(q)|^2.
\]
The first two terms are nonnegative, thus 
\begin{equation}\label{eq:k1-bound}
0\ge\sum_{5\nmid q}
 \left(qs\cot\frac{\pi q}{5}-\frac25\right)
 \left(qs\cot\frac{\pi q}{5}-1\right)
 |\widehat\alpha(q)|^2.
\end{equation}

Next we estimate the coefficients in \eqref{eq:k1-bound}.
Since $q\cot(\pi q/5)$ is even, it is sufficient to consider
$q>0$. 

For $q=1$, the coefficient is greater than $-1/10$ by the inequality  $\left(x-\frac25\right)(x-1)>-\frac1{10}$, which holds for all real $x$.
 
For $q\equiv3,4\pmod5$, the multiplier
$qs\cot(\pi q/5)$ is negative, so the coefficient
is greater than $2/5$.

For $q=2$, Lemma~\ref{lem:k1-angles} gives
\[
 2s\cot\frac{2\pi}{5}
 \le\frac{3-\sqrt5}{2},
\]
so the coefficient is at least $(\sqrt5-2)^2/5$.

For $q\ge6$ with $q\equiv1,2\pmod5$,  by  Lemma~\ref{lem:k1-angles} we have
\[
 qs\cot\frac{\pi q}{5}
 \ge 7s\cot\frac{2\pi}{5}
 >\frac{63}{20}\cot\frac{2\pi}{5},
\]
and thus the coefficient is at least $(\sqrt5-2)^2/5$.

Therefore, \eqref{eq:k1-bound} implies
\begin{equation}\label{eq:alpha-up-bound}
0\geq \frac{(\sqrt5-2)^2}{5}\sum_{q\equiv\pm2\pmod5}|\widehat\alpha(q)|^2
-\frac1{10}\left(|\widehat\alpha(1)|^2+
|\widehat\alpha(-1)|^2\right).
\end{equation}
By Parseval's identity and
the reality of $\alpha$,
\[
 \sum_{q\equiv\pm2\;(\mathrm{mod}\,5)}
 |\widehat\alpha(q)|^2=2\|z_2\|_2^2,
 \qquad
 |\widehat\alpha(1)|^2+|\widehat\alpha(-1)|^2
 \le2\|z_1\|_2^2,
\]
and thus \eqref{eq:alpha-up-bound} gives
\begin{equation}\label{eq:k1-dynamic-bound}
 2(\sqrt5-2)^2 \|z_2\|_2^2 \le \|z_1\|_2^2.
\end{equation}
But the geometric bound
\eqref{eq:k1-geometric-bound} gives
\[
 \|z_1\|_2^2
 \le(\sqrt5-2)^2\|z_2\|_2^2.
\]
Therefore, $z_1=z_2=0$ and Lemma~\ref{lem:constant-angle} completes the proof.
\end{proof}


\begin{thebibliography}{99}

\bibitem{ARSY}
M.~A.~Alfonseca, D.~Ryabogin, A.~Stancu, and V.~Yaskin,
\emph{On flotation, stability and related questions: A survey}, New Probes into Discrete and Convex Geometry, 2026. To appear.


\bibitem{AKZ16}
O.~Asipchuk, M.~Kosmakov, and P.~Zatitskii,
\emph{Rigidity in the planar Ulam floating body problem with perimetral density
$\sigma=1/6$}, Int. Math. Res. Not., 2026 (18), rnag203.

\bibitem{AKZ18}
O.~Asipchuk, M.~Kosmakov, and P.~Zatitskii,
\emph{Rigidity in the planar Ulam floating body problem with perimetral densities
$\sigma=1/8,3/8$ under central symmetry}, arXiv:2607.08338, 2026.

\bibitem{Auerbach}
H.~Auerbach,
\emph{Sur un probl\`eme de M.~Ulam concernant l'\'equilibre des corps flottants},
Studia Math. \textbf{7} (1938), 121--142.

\bibitem{BMO1}
J.~Bracho, L.~Montejano, and D.~Oliveros,
\emph{A classification theorem for Zindler carrousels},
J. Dynam. Control Systems \textbf{7} (2001), no.~3, 367--384.

\bibitem{BMO2}
J.~Bracho, L.~Montejano, and D.~Oliveros,
\emph{Carousels, Zindler curves and the floating body problem},
Period. Math. Hungar. \textbf{49} (2004), no.~2, 9--23.

\bibitem{MR2026}
Y.~Martinez-Maure and  D.~Rochera,
\emph{Fractal solutions to the two-dimensional floating body problem}, arXiv:2611.11309, 2026.

\bibitem{ScottishBook}
R.~D. Mauldin, editor,
\emph{The Scottish Book: Mathematics from the Scottish Caf\'e, with Selected
Problems from the New Scottish Book}, 2nd updated and enlarged ed.,
Birkh\"auser/Springer, Cham, 2015.

\bibitem{Ryabogin2022}
D.~Ryabogin,
\emph{A negative answer to Ulam's Problem 19 from the Scottish Book},
Ann. of Math. (2) \textbf{195} (2022), no.~3, 1111--1150.

\bibitem{Salkowski1934}
E.~Salkowski,
\emph{Eine kennzeichnende Eigenschaft des Kreises},
Sitzungsber. Heidelberger Akad. Wiss. Math.-Natur. Kl.  (1934), 57--62.


\bibitem{Wegner2003}
F.~J.~Wegner,
\emph{Floating bodies of equilibrium},
Stud. Appl. Math. \textbf{111} (2003), 167--183.

\bibitem{Wegner2007}
F.~J.~Wegner,
\emph{Floating bodies of equilibrium in 2D, the tire track problem and electrons
in a parabolic magnetic field}, arXiv:physics/0701241, 2007.

\bibitem{Zindler}
K.~Zindler,
\emph{\"Uber konvexe Gebilde. II. Teil},
Monatsh. Math. Phys. \textbf{31} (1921), 25--56.

\end{thebibliography}
\end{document}